\documentclass{article}
\usepackage{amssymb, amsmath, amsthm, mathtools, graphicx, standalone,enumerate, tikz}
\usepackage[left=1in,top=1in,right=1in]{geometry}
\usepackage{hyperref}
\usepackage[
style=alphabetic,
sorting=ynt,
defernumbers=true
]{biblatex}
\usepackage{xcolor}
\usepackage{soul}
\setstcolor{blue}
\usepackage{authblk}
\DeclareBibliographyCategory{cited}
\AtEveryCitekey{\addtocategory{cited}{\thefield{entrykey}}}
\nocite{*}

\theoremstyle{plain}
      \newtheorem{theorem}{Theorem}[section]
      \newtheorem{lemma}[theorem]{Lemma}

      \newtheorem{corollary}[theorem]{Corollary}
      \newtheorem{proposition}[theorem]{Proposition}
      
\theoremstyle{plain}
      \newtheorem{definition}[theorem]{Definition}
\theoremstyle{remark}

\newcommand{\Aa}{\mathcal{A}}
\newcommand{\Bb}{\mathcal{B}}

\newcommand{\Ll}{\mathcal{L}}

\newcommand{\Hh}{\mathcal{H}}
\renewcommand{\Rn}{\mathbb{R}^n}
\newcommand{\RR}{\mathbb{R}}
\newcommand{\cl}{c\ell}

\newcommand{\join}{\vee}

\newcommand{\Rec}{\operatorname{Rec}}

\title{Region level via centralization for hyperplane arrangements and beyond}

\author{Finn Southerland \thanks{University of California, San Diego. Email: \texttt{fsoutherland@ucsd.edu}.} \and Lani Southern \thanks{University of California, San Diego. Email: \texttt{lsouthern@ucsd.edu}.} \and Su Zhou \thanks{University of California, San Diego. Email: \texttt{suzhou@ucsd.edu}.}}

\date{}

\begin{document}

\maketitle

\begin{abstract}
In \cite{Zaslavsky_2003}, Zaslavsky showed how to compute the number $r_\ell(\Aa)$ of regions of a real hyperplane arrangement $\Aa$ with a given \emph{level}, refining his well known enumeration of regions and relatively bounded regions \cite{zaslavsky_1975}. We restate this theorem in terms of a construction called the \emph{centralization} of $\Aa$, give a bijective proof, and then apply it in two ways to answer active questions concerning the concept of level.
Firstly, a consequence of this enumeration is that $r_\ell(\Aa)$ depends only on the intersection poset $\Ll(\Aa)$, such that both $r_\ell$ and centralization can be defined in the more general setting of geometric semilattices. In this context we derive a very general expression for the characteristic polynomial of a geometric semilattice with several interesting corollaries. Secondly, recent investigations into the phenomenon of level have made little use of Zaslavsky's level-counting theorem, but it can be applied to obtain or generalize many of their results. In particular we show how exponential generating function identities \cite{chen2025level} \cite{chen_2026} and an expression giving the characteristic polynomial in terms of $r_\ell$ \cite{zhang_2024} can be derived for deformations of the braid arrangement.
\end{abstract}

\section{Introduction}
A \emph{hyperplane arrangement} $\Aa$ in $\Rn$ is a finite collection of affine hyperplanes.
The complement of $\Aa$ consists of connected components called \emph{regions} separated by the hyperplanes of $\Aa$. In a celebrated result (Theorem~\ref{Zaslavsky}), Zaslavsky showed how to compute the numbers $r(\Aa)$ of regions and $b(\Aa)$ of relatively bounded regions for a given arrangement via its characteristic polynomial $\chi_\Aa(t)$ \cite{zaslavsky_1975}. In this paper we are concerned with a refinement of this enumeration where regions are partitioned by the \emph{level} statistic, which captures the degree to which a region is unbounded. Ehrenborg first defined level in \cite{ehrenborg_level}\footnote{Though this paper was not published until 2019, an earlier version existed by 2001 when it was cited in \cite{Zaslavsky_2003}.} and enumerated faces of the extended Shi arrangements by both dimension and level. Subsequently in \cite{Zaslavsky_2003}, Zaslavsky gave a formula for the number of faces of $\Aa$ of dimension $d$ and level $\ell$ in terms of the characteristic and Whitney polynomials of some related arrangements, though he used the equivalent (up to a shift) concept of \emph{ideal dimension}. Independently, in \cite{Armstrong_Rhoades_2011} Armstrong and Rhoades gave another equivalent definition, which they called \emph{degrees of freedom}. We will use the term level along with this last definition.
\begin{definition}
    The \emph{level} of a region $R$ is $\ell(R) = \dim \Rec(R)$ where $\Rec(R)$ is the \emph{recession cone} of $R$:
\[\Rec(R) = \{ v \in \Rn : R + v \subseteq R\}.\]
\end{definition}

The reader may wish to check that $\Rec(R)$ is indeed closed under nonnegative linear combination, and is therefore a cone in the geometric sense.
Alternatively, the level of $R$ is the smallest dimension of a subspace such that the distance from a point of $R$ to the subspace is bounded.
A region has recession cone $\{0\}$ if and only if it is bounded, so bounded regions are those where $\ell(R) = 0$. A region with $\ell(R) = 1$ looks like a (one or two-way infinite) tube, a region with $\ell(R) = 2$ looks like some part of a ``slab", and so on. We write $r_\ell(\Aa)$ for the number of regions of $\Aa$ of level $\ell$.

\begin{figure}
    \centering
    \includegraphics[width=0.2\linewidth]{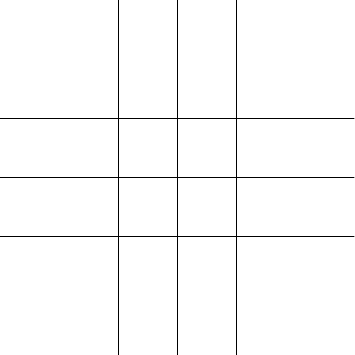}
    \hspace{1cm}
    \includegraphics[width=0.2\linewidth, trim={0, .5cm, 0, 0},clip]{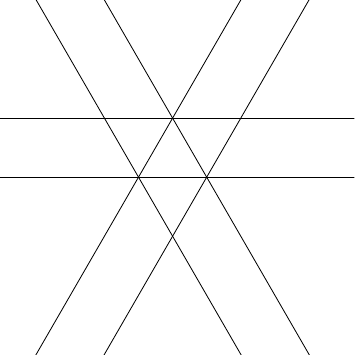}
    \caption{A pair of arrangements with the same characteristic polynomial, but different values of $r_1$ and $r_2$.}
    \label{fig:level_ex}
\end{figure}

A consequence of Zaslavsky's region counting theorem is that for any pair of hyperplane arrangements $\Aa_1,\Aa_2$ with the same characteristic polynomial, we also have $r(\Aa_1) = r(\Aa_2)$ and $b(\Aa_1) = b(\Aa_2)$. Although the terminology is not common, we might therefore call these numbers \emph{characteristic invariants}. A larger class of invariants is the \emph{combinatorial invariants}, which are determined by the intersection poset $\Ll(\Aa)$. The pair of arrangements in Fig. \ref{fig:level_ex} show that $r_\ell$ is not a characteristic invariant, but we will see that it is a combinatorial invariant as a consequence of Theorem \ref{degrees of freedom theorem} \cite[Theorem 3.1]{Zaslavsky_2003}. Zaslavsky used different (though closely related) definitions; we express $r_\ell(\Aa)$ as a function of $\Ll(\Aa)$ using the \emph{centralization} of $\Aa$, to our knowledge first defined and shown to be a combinatorial invariant itself in \cite{Budur_2022}.

\begin{definition}
    For a hyperplane $H = \{ x \in \Rn:w\cdot  x = a\}$, the \emph{centralization} of $H$ is \[\underline{H} := \{ x \in \Rn : w \cdot  x = 0\},\] the unique parallel hyperplane passing through the origin. For an arrangement $\Aa = \{H_1,\ldots, H_k\}$ the \emph{centralization} of $\Aa$ is $\underline{\Aa} := \{\underline{H_1}, \ldots, \underline{H_k}\}$.
\end{definition}

As suggested by its name, $\underline{\Aa}$ is a \emph{central} arrangement, since the intersection of all its hyperplanes contain the origin. The concept of centralization is implicit in many places in the literature (as we will see in Section \ref{sec:combinatorics}), but has received relatively little attention. Because of the simplicity of this definition and its naturality in the context of level, centralization may be worthy of further study. With this tool in hand, as well as the \emph{restriction} $\Aa^V$ and \emph{localization} $\Aa_V$ of an arrangement to a subspace $V$ (see Section \ref{sec:background} for definitions), we can state a version of Zaslavsky's level-counting result.

\begin{theorem} \cite[Equation 3.1]{Zaslavsky_2003}
\label{degrees of freedom theorem}
Let $\Aa$ be a hyperplane arrangement in $\Rn$, and $0 \leq \ell \leq n$. Then the number of regions of $\Aa$ of level $\ell$ is a combinatorial invariant\footnote{Zaslavsky does not consider combinatorial invariance in \cite{Zaslavsky_2003}, but it is a direct consequence of the formula. This seems not to have been appreciated as the question was thought to be open by some much more recently and provided the impetus for our study.}. In particular, 
\begin{equation} \label{eqn:main}
        r_\ell(\Aa) =
\sum_{V\in \Ll_{n-\ell}(\underline{\Aa})} r(\underline{\Aa}^{V}) b(\Aa_V).
    \end{equation}
\end{theorem}
\begin{proof}
    The correctness of Equation \eqref{eqn:main} will be established in Theorem \ref{correctness}. 
    By Zaslavsky's Theorem (\ref{Zaslavsky}), each summand can be calculated from $\Ll(\underline{\Aa}^V)$ and $\Ll(\Aa_V)$. Theorem \ref{combinatorial_inv} says that these, as well as $\Ll(\underline{\Aa})$, are combinatorial invariants, completing the proof.
\end{proof}
We remark that for an essential arrangement in the case $\ell = 0$, Equation \eqref{eqn:main} reduces to simply $b(\Aa)$, as $\Ll_n(\underline{\Aa})=\{0\}$, $\underline{\Aa}^0$ is an empty arrangement with 1 region, and $\Aa_0 = \Aa$. In the other extreme, the case $\ell=n$ reduces to $r_n(\Aa) = r(\underline{A})$, which reflects the general facts that $r_n(\Aa)$ is constant under translation of the hyperplanes of $\Aa$, and that every region of a central arrangement has level $n$.

One way our statement of Theorem \ref{degrees of freedom theorem} differs from Zaslavsky's in \cite{Zaslavsky_2003} is that he gives a more general enumeration of faces with dimension $d$ and level $\ell$. This enumeration is obtained from Theorem \ref{degrees of freedom theorem} by summing over all restrictions $\Aa^F$, for $F \in \Ll(\Aa)$. While in this paper we do not otherwise consider counting faces of dimension other than $n$, we give this extension to help clarify the relationship between our notation and concepts and those used by Zaslavsky.
\begin{corollary}
    Let $\Aa$ be a hyperplane arrangement in $\RR^n$. For any $0 \leq \ell \leq d \leq n$, the number of faces of $\Aa$ of dimension $d$ and level $\ell$ is
\begin{equation}
    f_{d,\ell}(\Aa) = \sum_{F \in \Ll_{n-d}(\Aa)}\sum_{V \in \Ll_{d - \ell}(\underline{\Aa}^F)} r(\underline{\Aa}^V)b(\Aa^F_V).
\end{equation}
\end{corollary}
\begin{proof}
    As every face of dimension $d$ is contained in a unique flat $F \in \Ll(\Aa)$ of the same dimension, it suffices to sum Equation \ref{eqn:main} over all such flats, replacing $\Aa$ by $\Aa^F$. In each summand, it is easy to check that restriction commutes with both centralization and localization, and that $(\underline{\Aa}^F)^V = \underline{\Aa}^V$ for each $V \in \Ll_{d - \ell}(\underline{\Aa}^F)$.
\end{proof}

The statement of this corollary still differs from Zaslavsky's in that he gives it in terms of a bivariate generating function, which allows collapsing the outer sum into the definition of the Whitney polynomial of $\Aa$. More significantly, Zaslavsky uses the concepts of ideal dimension and projectivization, which are equivalent but not identical to our level and centralization. To compare ideal dimension with level, the reader should note that all sets that we consider the level of are relatively open polyhedra, and so have level one greater than their ideal dimension \cite[Proposition 7.4]{Zaslavsky_2003}. To compare the other concepts, we have the isomorphisms of posets $\Ll(\Hh^\infty) \simeq \Ll(\underline{\Aa})$, $\Ll(\Hh^s) \simeq \Ll(\underline{\Aa}^V)$, and $\Ll(\Hh(s)) \simeq \Ll(\Aa_V)$ in Zaslavsky's Equation 1.3 and our Equation \eqref{eqn:main} respectively.

Zaslavsky's proof is not explicitly bijective, though it could be modified to be so. To improve this situation and to show how centralization can be used in place of projectivization, we give a bijective proof for Theorem \ref{degrees of freedom theorem} between Sections \ref{sec:lemmas} and \ref{sec:main}. In the first, we develop three geometric lemmas regarding recession cones. In the second, we give a bijection for each summand of Equation \eqref{eqn:main}, stated as Theorem \ref{bijection}. This bijection is also useful later, in obtaining Theorem \ref{chi_thm}.

In \ref{sec:combinatorics}, we prove the combinatorial invariance of $r_\ell(\Aa)$ by giving an explicit algorithm for computing $r_\ell(\Aa)$ from $\Ll(\Aa)$. In particular, we show how to construct $\Ll(\underline{\Aa})$ given $\Ll(\Aa)$, and how to evaluate Equation \eqref{eqn:main} from the two posets. Though the poset $\Ll(\Aa)$ is an example of a geometric semilattice, not every geometric semilattice is the intersection lattice of a real arrangement $\Aa$. Nonetheless we can apply the algorithm represented by Equation \eqref{eqn:main} to any geometric semilattice $M$ and obtain the geometric lattice $\underline{M}$, which we call the centralization of $M$, and a distribution $r_\ell(M)$. Finding a combinatorial, geometric, or algebraic interpretation of this distribution for an abstract semilattice remains an open question, since there is no concept of a ``region'' of a geometric semilattice.
However, we are able to prove results about the distribution at this level of generality, and exploring its properties is the focus of Section \ref{sec:semilattice}. Our main result is an expression for the characteristic polynomial of a geometric semilattice.

\begin{theorem} \label{chi_semilattice}
    For a geometric semilattice $M$, we have 
    \begin{equation}
        \chi_M(t) = \sum\limits_{S \in \underline{M}} \chi_{\underline{M}^S}(t) \chi_{M_S}(1).
    \end{equation}
\end{theorem}

This theorem has several interesting consequences. For one, it gives us the identity $r(M) = \sum_{\ell = 0}^n r_\ell(M)$ for geometric semilattices, where the lack of a known combinatorial interpretation for these numbers makes this a nontrivial fact. When $M = \Ll(\Aa)$ for an arrangement $\Aa$ we can rearrange this theorem into a sum over the regions of $\Aa$, and this is Theorem \ref{chi_thm}. Finally, for specific cases where the geometric lattice $\underline{M}$ has the property that its principal order filters $\underline{M}^S$ depend only on $\operatorname{rank}(S)$, terms of the sum can be collected by rank, giving Corollary \ref{cor:uniform}. A motivating example of this phenomenon is the \emph{braid arrangement}, whose deformations have been studied extensively in the context of level, and whose intersection poset satisfies this property. The reader should keep in mind that geometric lattices and semilattices are conceptually equivalent to (simple) matroids and semimatroids - therefore centralization and the level distribution are also defined in that setting.

Other results (e.g. \cite{chen_2026}, \cite{ehrenborg_level}), on counting regions by level have also been concerned with various deformations of the braid arrangement. In Section \ref{sec:exp} we apply Theorem \ref{degrees of freedom theorem} to the class of exponential arrangements, generalizing some of these. Here Corollary \ref{cor:uniform} recovers the results of \cite{zhang_2024}, which in turn generalize results in \cite{chen_2026} and \cite{ehrenborg_level}, expressing the characteristic polynomial in terms of the region level counts for many arrangements.

\section{Background}\label{sec:background}
\subsection{Real hyperplane arrangements}
Here we establish the rest of the definitions we will need to prove the upcoming results. A more thorough treatment of most of these can be found in \cite{stanley2007introduction} or \cite{orlik_terao_1992}.
Throughout this paper all hyperplane arrangements should be assumed to be in $\RR^n$. We order the hyperplanes of an arrangement $\Aa = \{H_1,\ldots,H_k\}$ and assume that some $w_i \in \Rn$ and $a_i \in \RR$ have been selected so that $H_i = \{x \in \Rn:w_i \cdot x = a_i\}$. For each $H_i \in \Aa$ we have a decomposition $\Rn = H_i \sqcup H_i^+ \sqcup H_i^-$, where $H_i^+$ and $H_i^-$ are the open half spaces where $w_i \cdot x > a_i$ and $w_i \cdot x < a_i$, respectively. A \emph{region} of $\Aa$ is a a nonempty intersection $R = \cap_{i=1}^k H_i^{\sigma(R)_i}$, where $\sigma(R) \in \{\pm1\}^k$ is the \emph{sign vector} of the region. A \emph{face} of $\Aa$ is a similar intersection where the $i$th term is drawn from $\{H_i,H_i^+,H_i^-\}$.

The \emph{intersection poset} $\Ll(\Aa)$ is the set of nonempty intersections of the hyperplanes of $\Aa$, which are called \emph{flats} and ordered by reverse inclusion. The poset is a geometric (meet) semilattice \cite{ardila2004_semimat} with minimum element $\hat 0 = \Rn$ being the intersection over the empty set, and atoms $H \in \Aa$. When the intersection of all faces of $\Aa$ is nonempty we call $\Aa$ \emph{central} and observe that $\Ll(\Aa)$ is a lattice whose maximum element we denote by $\hat 1$. The rank function of $\Ll(\Aa)$ is $\operatorname{rank}(F) = \operatorname{codim}(F)$. By $\Ll_r(\Aa)$ we denote the elements of $\Ll(\Aa)$ with rank $r$.

The \emph{characteristic polynomial} of $\Aa$ is \[\chi_\Aa(t) = \sum\limits_{F \in \Ll(\Aa)} \mu(\hat 0, F)t^{\dim(F)},\] where $\mu$ is the M\"obius function on $\Ll(\Aa)$. We may have $\operatorname{rank}(\Ll(\Aa)) < n$, and in this case we call $\Aa$ \emph{nonessential}. Letting $W \subseteq \Rn$ be the span of the normal vectors $w_i$, the arrangement obtained by restricting $\Aa$ to $W$ is called the \emph{essentialization} of $\Aa$, and denoted $\operatorname{ess}(\Aa)$. This arrangement is \emph{essential}, that is rank(ess$(\Aa))=\dim(W)$, and combinatorially equivalent to $\Aa$, that is $\Ll(\operatorname{ess}(\Aa)) \simeq \Ll(\Aa)$. There is a natural bijection $R \leftrightarrow R \cap W$ between regions of $\Aa$ and $\operatorname{ess}(\Aa)$, and regions of $\Aa$ corresponding to bounded regions of $\operatorname{ess}(\Aa)$ are called \emph{relatively bounded}.

\begin{theorem}{(Zaslavsky \cite{zaslavsky_1975})}
\label{Zaslavsky}
 For a hyperplane arrangement $\Aa$ in $\Rn$, let $r(\Aa)$ be the number of regions of $\Aa$, and $b(\Aa)$ the number of relatively bounded regions. Then
 \[
r(\Aa) = (-1)^n\chi_\Aa(-1), \hspace{10pt}
b(\Aa) = (-1)^{\operatorname{rank}(\Aa)}\chi_\Aa(1).
\]
\end{theorem}

We make use of two notions of subarrangements of $\Aa$, both corresponding to certain subposets of $\Ll(\Aa)$. 

\begin{definition}
    The \emph{restriction} of an arrangement $\Aa$ to an affine subspace $V$ (often a flat of $\Aa$) is \[\Aa^V := \{H \cap V : H \in \Aa\} \setminus \{\emptyset, V\}.\]
\end{definition}

When $V$ is a flat of $\Aa$, $\Ll(\Aa^V)$ is simply the principal order filter $\Ll(\Aa)^V = \{W \in \Ll(\Aa) : V \leq W\}$.

\begin{definition} \label{localization}
    For a linear subspace $V \subseteq \Rn$, the localization of $\Aa$ at $V$ is \[\Aa_V := \{H \in \Aa : V \subseteq H \text{ or } V \cap H = \emptyset\}.\] 
\end{definition}

This is the order ideal of $\Ll(\Aa)$ consisting of flats parallel to $V$. This notation is made in analogy to the principal order ideal $P_s = \{t \in P: t \leq s\}$ for a poset $P$ and element $s \in P$.

\subsection{Geometric lattices and semilattices}

In proving combinatorial invariance of the level distribution, we will identify $\Ll(\Aa)$ and $\Ll(\underline{\Aa})$ as a geometric semilattice and a geometric lattice, respectively. 

\begin{definition}
A ranked lattice $L$ is \emph{semimodular} if for every $s,t \in L$ we have 
\[\operatorname{rank}(s) + \operatorname{rank}(t) \geq \operatorname{rank}(s \wedge t) + \operatorname{rank}(s \vee t).\]
If for all $s \in L$ we have that $s = \bigvee\limits_{a \in S} a$ for some $S \subseteq L_1$, that is $s$ is the join of some subset of the atoms of $L$, we say $L$ is \emph{atomistic}. If $L$ is both semimodular and atomistic we say it is a \emph{geometric lattice}.
\end{definition}

\begin{definition}
\label{def:geometric semilattice}
    Let $M$ be a ranked meet semilattice, then we say that a set $S$ of atoms of $M$ is \emph{independent} if they have an upper bound $\bigvee S$ and if $r(\bigvee S) = |S|$. We call $M$ a \emph{geometric semilattice} if it satisfies
    \begin{enumerate}[(a)]
        \item every principal order ideal of $M$ is a geometric lattice, and
        \item whenever $S$ is an independent set of atoms of $M$ and $r(t) < r(\vee S)$, there is some $a \in S$ such that $a \not\leq t$ and $t \vee a$ exists.
    \end{enumerate}
\end{definition}

The intersection posets of hyperplane arrangements are examples of geometric semilattices, and intersection posets of central arrangements are geometric lattices. Note that geometric lattices and geometric semilattices correspond to simple matroids and semimatroids respectively \cite{ardila2004_semimat}, though we do not use that language in this paper. More discussion of geometric lattices and semilattices can be found in \cite{Wachs1985}.

\section{Recession cones} \label{sec:lemmas}
In this section we establish three lemmas about recession cones that we will use in the next to establish the bijection at the heart of Theorem \ref{degrees of freedom theorem}.
The first two of these are alternative characterizations of $\Rec(R)$, one that is more explicit than the definition, and one that allows us to focus on a single point of $R$.
Our first characterization is in terms of the sign vector $\sigma(R)$, which we observe satisfies $\sigma(R)_i(w_i \cdot x - a_i) > 0$ for all $x \in R$. Note that $\sigma(R)$ depends on the choices of $w_i$s, and when considering a single region we will often assume that these have been made so that $\sigma(R)_i = +1$ for all $1 \leq i \leq k.$

\begin{lemma} \label{rec_char_1}
For an affine arrangement $\Aa = \{H_1,\ldots,H_k\}$ and fixed region $R$ of $\Aa$, we have 
\[\Rec(R) = \{v \in \Rn : \text{for all } 1 \leq i \leq k, \sigma(R)_i(w_i \cdot v) \geq 0\}.\]
\end{lemma}
\begin{proof}
    Assume that $\sigma(R)_i = 1$ for all $i$. Then for all $x \in R$ and $v \in \Rn$ such that $w_i \cdot v \geq 0$, we have \[w_i \cdot (x + v) = w_i\cdot x + w_i \cdot v \geq w_i \cdot x > a_i,\] so $x + v \in R$. On the other hand, suppose that for some $v \in \Rec(R)$ there is an $i$ such that $w_i \cdot v < 0$. For an arbitrary $x \in R$, let $c = w_i \cdot x - a_i$. Since $\Rec(R)$ is a cone, we may scale $v$ to obtain some $v' \in \Rec(R)$ such that $w_i \cdot v' < -c$. Then we have \[w_i \cdot (x + v')  = w_i \cdot x + w_i \cdot v'  < (c + a_i) - c= a_i,\] so $x + v' \notin R$, a contradiction.
\end{proof}

This characterization illustrates how $\underline{\Aa}$ becomes involved when we study recession cones. Where $R$ is an intersection of open half-spaces defined by the $H_i$, the recession cone $\Rec(R)$ is the corresponding intersection of closed half-spaces defined by the $\underline{H_i}$. Such an intersection is a \emph{face} of $\underline{\Aa}$, i.e. a region of a restriction of $\underline{\Aa}$. This will allow us to break up the set of regions of $\Aa$ by considering those regions whose recession cone is a particular face of $\underline{\Aa}$, of which there may be multiple.

\begin{lemma} \label{rec_char_2}
    For any $x \in R$, we have \[\Rec(R) = \{v \in \Rn : \text{for all } c \in \RR^+, x + cv \in R\}.\]
\end{lemma}
\begin{proof}
    Call the set on the right-hand side $X$, then the containment $\Rec(R) \subseteq X$ follows from the fact that $\Rec(R)$ is a cone. On the other hand, letting $v \in X$ and $y \in R$, it suffices to show that $y + v \in R$. Since $R$ is open, we can pick $y' \in R$ such that $y$ lies on the interval strictly between $x$ and $y'$. Then the ray from $y'$ through $y + v$ and the ray from $x$ through $x + v$ intersect at a point $x + cv$ for some $c \in \RR^+$, so this intersection is in $R$ by the definition of $X$. By convexity of $R$, $y + v \in R$ and so we have $v \in \Rec(R)$.
\end{proof}

Finally, it is useful to pick out an element of $\Rec(R)$ for which as many of the inequalities of Lemma \ref{rec_char_1} are strict as is possible.

\begin{lemma} \label{rec_interior}
    For some region $R$ of $\Aa$, let $V = \operatorname{span}(\Rec(R))$. Then there is some $v \in \Rec(R)$ such that for every $i$ where $\underline{H_i}^V$ is a hyperplane of $\underline{\Aa}^V$, we have $\sigma(R)_i(w_i \cdot v) > 0$.
\end{lemma}
\begin{proof}
    Again, assume $\sigma_i(R) = 1$ for all $i$. The condition given is equivalent to $V \not\subseteq \underline{H_i}$, so for each such $i$ we may choose $v_i \in \Rec(R) \setminus \underline{H_i}$. Then $w_i \cdot v_i > 0$, and $w_j \cdot v_i \geq 0$ for all $1 \leq j \leq k$. Letting $v = \sum v_i$ we observe that $w_i \cdot v > 0$ for all $i$ satisfying the condition, and $v \in \Rec(R)$ by the fact that it is a cone.
\end{proof}

Geometrically this lemma says that the hyperplanes $H_i$ such that $\underline{H_i}^V \in \underline{\Aa}^V$ are exactly those that a point in $R$ can get arbitrarily far away from. In practice we will use this lemma with scaling to select $v$ so that $w_i \cdot v > a$ for some constant $a \in \mathbb{R}^+$. 

\section{Enumerating regions of level \texorpdfstring{$\ell$}{l}} \label{sec:main}
Our intuition behind the equation of Theorem \ref{degrees of freedom theorem} involves ``zooming out" on $\Aa$ so that its hyperplanes appear arbitrarily close to linear. This transforms $\Aa$ into $\underline{\Aa}$, and (the closure of) each region of $\Aa$ into its recession cone, which as we saw in the last section is a face of $\underline{\Aa}$. Taking the span of the recession cone gives a flat of $\underline{\Aa}$, and so to each region $R$ we associate the flat $V_R \in \Ll(\underline{\Aa})$, where $\ell(R) = \dim(V_R)$. Equation \eqref{eqn:main} simply counts regions of level $\ell$ by summing over flats of the correct dimension, and it remains to determine the number of regions $R$ such that $V_R = V$. We will show that this is $r(\underline{\Aa}^V)b(\Aa_V)$, where the factor of $r(\underline{\Aa}^V)$ comes from distinguishing the particular face $\Rec(R) = F\subseteq V$, and the central claim of the equation is that $b(\Aa_V)$ is the number of regions such that $\Rec(R) = F$ for any face $F$ of $\underline{\Aa}$. To prove this, we will give a map
\begin{equation}
    \Phi_V: \left\{\text{Regions of } \Aa \text{ such that } V_R = V\right\} \to \{\text{Regions of } \underline{\Aa}^V\} \times \{\text{Relatively bounded regions of } \Aa_V\}
\end{equation}
for each $V \in \Ll(\underline{\Aa})$ and provide its inverse $\Psi_V$.

We define $\Phi_V$ using the sign vector of $R$. The arrangements $\underline{\Aa}^V$ and $\Aa_V$ each include hyperplanes corresponding to a subset of the hyperplanes of $\Aa$; in the case of $\underline{\Aa}^V$ these are the $\underline{H_i}^ V$s (which are no longer unique). This allows us to define the region $R_u$ of $\underline{\Aa}^V$ and the relatively bounded region $R_b$ of $\Aa_V$ as the regions of their respective arrangements with the ``same" sign vectors as $R$, that is:
\begin{equation}
\Phi_V : R \mapsto \begin{cases}
R_u = \{x \in V : \sigma(R)_i(w_i \cdot x) > 0 \text{ for each $i$ such that } \underline{H_i}^V \in \underline{\Aa}^V\} \\ 
R_b = \{x \in \Rn : \sigma(R)_i(w_i \cdot x) > a_i\text{ for each $i$ such that }H_i \in \Aa_V \}
\end{cases}.
\end{equation}
Note that from their definitions, $R_u$ and $R_b$ are either empty or regions of their respective arrangements. We will show that these are in fact regions and that $R_b$ is relatively bounded.

Given $R_u,$ $R_b$, we claim that the inverse map $\Psi_V$ is given by combining the sign vectors of $R_u$ and $R_b$ into a single sign vector and letting $R$ be the corresponding region. That is,
\begin{equation}
    \Psi_V : (R_u,R_b) \mapsto R = \left\{x \in \Rn : \begin{cases}
        \sigma(R_u)(w_i \cdot x - a_i) > 0 & \text{for $i$ such that } \underline{H_i}^V \in \underline{\Aa}^V \\
        \sigma(R_b)(w_i \cdot x_i - a_i) > 0 & \text{for $i$ such that }H_i \in \Aa_V
    \end{cases}\right\}.
\end{equation}
Since $\underline{H_i}$ either intersects non-trivially with or contains $V$, exactly one of $\underline{\Aa}^V$ and $\Aa_V$ has a hyperplane corresponding to $H_i$, and each $H_i$ appears exactly once in the conditions given. We will show that $R$ is a region of $\Aa$ (i.e. nonempty) and that $V_R = V$. Provided that these maps are well-defined, it will not be hard to show that they are inverses, since they essentially operate by splitting and recombining sign vectors.

\begin{theorem} \label{bijection}
For a hyperplane arrangement $\Aa$ and a given $V \in \Ll(\underline{\Aa})$, the number of regions $R$ of $\Aa$ such that $V_R = V$ is $r(\underline{\Aa}^V)b(\Aa_V)$. In particular, $\Phi_V$ and $\Psi_V$ are (well-defined) mutually inverse bijections.
\end{theorem}

\begin{proof}
Consider a particular region $R$ such that $V_R = V$. Throughout, we will assume that the $w_i$ have been chosen so that $\sigma(R)_i = +1$ for all $i$ and additionally that the hyperplanes of $\Aa$ are labeled such that $\underline{\Aa}^V = \{\underline{H_i}^V : 1 \leq i \leq k_u\}$ and $\Aa_V = \{H_i : k_u < i \leq k\}$ for some $k_u$.

We begin by showing that $\Phi_V$ is well-defined. Remember that showing $R_u,$ $R_b$ are regions only requires showing that they are nonempty. By Lemma \ref{rec_interior} we may pick a point $v_u \in \Rec(R)$ such that $w_i \cdot v_u > 1$ for all $1 \leq i \leq k_u$, and this is a point in $R_u$. Since $\Aa_V$ contains a subset of the hyperplanes of $\Aa$, we have $R \subseteq R_b$, so $R_b$ is also nonempty. It is therefore a region of $\Aa_V$, but we must also show that it is relatively bounded. All hyperplanes of $\Aa_V$ are parallel to $V$, so it suffices to show that $R_b \cap V^\perp$ is bounded. Suppose that $v_b \in \Rec(R_b \cap V^\perp)$, that is $w_i \cdot v_b \geq 0$ for all $k_u < i \leq k$. Note that $v_b \in V^\perp$. Now by scaling $v_b$ down we may assume that $w_i \cdot v_b > -1$ for all $i$ (in particular, for $1 \leq i \leq k_u$). Then adding the $v_u$ we chose earlier, we have $w_i \cdot (v_b + v_u) > 0$ for all $1 \leq i \leq k$, so $v_b + v_u \in \Rec(R)$. But now both $v_u$ and $v_u + v_b$ lie in $V$, so $v_b$ does as well. We have $v_b \in V \cap V^\perp = \{0\}$, and conclude that $\Rec(R_b \cap V^\perp) = \{0\}$ as desired.

Now we turn to $\Psi_V$. We continue to assume that our regions lie on the positive side of all hyperplanes, so we have $R = \{x \in \Rn : w_i \cdot x > a_i\}$. To show that this is a region of $\Aa$ (i.e. nonempty), let $x_b \in R_b$, and  $a > \max_{1 \leq i \leq k}(a_i - w_i \cdot x_b)$. Once again by Lemma \ref{rec_interior}, let $v_u \in V$ be such that $w_i \cdot v_u > a$ for all $1 \leq i \leq k_u$ and let $x = x_b  +v_u$. Then for $1 \leq i \leq k_u$, we have $w_i \cdot x = w_i \cdot x_b + w_i \cdot v_u > (w_i \cdot x_b) + (a_i - w_i \cdot x_b) = a_i$. For $k_u < i < k$, we have $w_i \cdot v_u = 0$, so $w_i \cdot x = w_i \cdot x_b > a_i$ as well. Therefore $x \in R$, so $R$ is a nonempty region of $\Aa$. We also need to show that $V_R = V$. In fact we will show that $\Rec(R) = \overline{R_u}$ (the closure of $R_u$), and since $R_u$ is a region of a central arrangement living in $V$ we have $\operatorname{span}(R_u) = V$.

For $x \in R_u$ we have $w_i \cdot x > 0$ for $1 \leq i \leq k_u$ by definition, and $w_i \cdot x = 0$ for $k_u < i \leq k$ as $R_u \subseteq V$. Taking the closure turns the strict inequalities non-strict, and comparing to Lemma \ref{rec_char_1} we see that $\overline{R_u} \subseteq \Rec(R)$. 

On the other hand, suppose $v \in \Rn \setminus \overline{R_u}$. To show that $v \notin \Rec(R)$, it is convenient to use Lemma \ref{rec_char_2} at $x$, selecting some $c>0$ such that $x + cv \notin R$. If $v \in V$, there must be some $1 \leq i \leq k_u$ such that $w_i \cdot v < 0$. Selecting $c>0$ so that $w_i \cdot cv < -(w_i \cdot x - a_i)$ we get $w_i \cdot (x + cv) < a_i$, so $x + cv \notin R$. Otherwise, we have $v \notin V$ and so it has nonzero projection $v'$ onto $V^\perp$. By the relative boundedness of $R_b$, we can let $c>0$ be such that $w_i \cdot (x_b + cv') - a_i < 0$ for some $k_u < i \leq k$. But $w_i \cdot v_u = 0$ and $w_i \cdot v = w_i \cdot v'$, so we have $w_i\cdot (x + cv) - a_i = w_i\cdot(x + cv) - a_i < 0$. Again, $x + cv \notin R$. This gives us $\Rec(R) \subseteq \underline{R_u}$ as desired.

Lastly, we must argue that $\Phi_V$ and $\Psi_V$ are inverses. To see this, note that going in either direction, $\sigma(R)_i = \sigma(R_u)_i$ for all $1 \leq i \leq k_u$ and $\sigma(R)_i = \sigma(R_b)_i$ for $k_u < i \leq k$. Since a region is uniquely identified by its sign vector, this establishes $\Psi_V(\Phi_V(R)) = R$ and $\Phi_V(\Psi_V(R_u,R_b))=(R_u,R_b)$.
\end{proof}

\begin{theorem} \label{correctness}
    Equation \eqref{eqn:main} counts the number of regions of $\Aa$ with level $\ell$.
\end{theorem}
\begin{proof}
    Lemma \ref{rec_char_1} shows that a region $R$ of $\Aa$ has level $\ell$ if and only if $\Rec(R)$ is a face of $\underline{\Aa}$ of dimension $\ell$, that is $\text{span}(\Rec(R)) = V_R$ for some $V_R \in \Ll(\underline{\Aa})_{n-\ell}$. Therefore \[r_\ell(\Aa) = \sum\limits_{V \in \Ll(\underline{\Aa})_{n-\ell}} \#\{\text{Regions of $\Aa$ such that $V_R = V$}\}.\] By Theorem \ref{bijection}, this summand is equal to $r(\underline{\Aa}^V)b(\Aa_V)$.
\end{proof}

\section{Centralization and cones} \label{sec:combinatorics}
In this section we will show that Equation \eqref{eqn:main} is a combinatorial invariant, and develop it as the definition of $r_\ell(M)$ when $M$ is a geometric semilattice. This amounts to showing that $\Ll(\underline{\Aa})$, $\Ll(\underline{\Aa}^V)$, and $\Ll(\Aa_V)$ can be constructed from only $\Ll(\Aa)$, without knowing the geometry of $\Aa$. First, we will examine the geometric construction of the \emph{cone} of $\Aa$, which contains the information we need and plays the role of \emph{projectivization} in \cite{Zaslavsky_2003}. Equivalent combinatorial constructions were given in \cite{Wachs1985} and \cite{ardila2004_semimat}, and \cite[Proposition 3.2]{Budur_2022} observes that the first of these shows how $\Ll(\underline{\Aa})$ is a combinatorial invariant. Arguments about equivalence and invariance have often been fairly inexplicit however, and here we aim to give a more explicit account. Both \cite{Wachs1985} and \cite{ardila2004_semimat} also give their constructions in terms of the matroid of $\Ll(\Aa)$, here we give one purely in the language of geometric (semi)lattices. On one hand, our goal is to make clear that \cite{Zaslavsky_2003} establishes the combinatorial invariance of $r_\ell(\Aa)$, since this seems not to have been well understood. Additionally, we use this construction to define the numbers $r_\ell(M)$ for an arbitrary geometric semilattice $M$ which gives an interesting generalization that we begin to study in the next section.

The \emph{cone} of $\Aa$ is a standard construction in the theory of arrangements (see, e.g. \cite[Section 1.1]{stanley2007introduction}).
\begin{definition}
For an arrangement $\Aa$ in $\Rn$ with hyperplanes $H_i = \{x \in \Rn : w_i \cdot x = b_i\}$, let \[cH_i = \{x \in \RR^{n+1} : x_1w_{i1} + \cdots + x_nw_{in} - x_{n+1}b_i = 0\},\] and let $H_0 = \{x \in \RR^{n+1} : x_{n+1} = 0\}$. Then the \emph{cone} of $\Aa$ is 
\[c\Aa = \{H_0, cH_1,\ldots,cH_k\}.\]
\end{definition}

One can imagine taking slices of $c\Aa$ parallel to $H_0$, perhaps starting with $\{x_{n+1} = 1\}$ and watching a movie that smoothly converts $\Aa$ into $\underline{\Aa}$. We compile some key observations about $c\Aa$ into the following lemma, which an interested reader should be able to verify without much difficulty.
\begin{lemma} \label{lem:cone_geo}
Let $\Aa$ be an arrangement in $\Rn$ and let $V \in \Ll(\underline{\Aa})$. The following isomorphisms of posets hold, with the right hand sides being induced subposets of $\Ll(c\Aa)$:
\begin{equation} \label{A_in_cone}
    \Ll(\Aa) \simeq \Ll(c\Aa) - \Ll(c\Aa)^{H_0},
\end{equation}
\begin{equation} \label{center_in_cone}
    \Ll(\underline{\Aa}) \simeq \Ll(c\Aa)^{H_0},
\end{equation}
\begin{equation} \label{res_in_cone}
    \Ll(\Aa_V) \simeq \Ll(c\Aa)_V - \Ll(c\Aa)^{H_0},
\end{equation}
where in Isomorphism \eqref{res_in_cone} we view $V$ as an element of $\Ll(c\Aa)$ through Isomorphism \eqref{center_in_cone}.
\end{lemma}

\begin{figure}
    \centering
    \includegraphics[width=8cm]{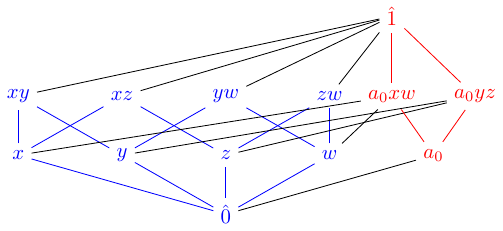}
    \caption{The cone of the lattice of flats of the matroid $U_{2,4}$, minus its maximum element. A copy of the original semilattice is in blue, with its centralization `above' it in red.}
    \label{fig:poset}
\end{figure}

This lemma shows that it will be sufficient to construct $\Ll(c\Aa)$ from $\Ll(\Aa)$ in such a way that $H_0$ is distinguished from the other atoms. Since $\Ll(\Aa)$ is a lattice when $\Aa$ is central, the cone $\Ll(c\Aa)$ is a geometric lattice. We will construct the cone of an arbitrary geometric semilattice, which will similarly be a geometric lattice. An example is shown in Figure~\ref{fig:poset}.

\begin{definition} \label{combo_cone}
    Let $M$ be a geometric semilattice with atoms $A = \{a_1,\ldots,a_k\}$. We construct a geometric lattice $cM$ which we call the \emph{cone} of $M$ with atoms $\{\{a_0\}, \{a_1\},\ldots,\{a_k\}\}$ for a distinguished additional `atom` $a_0$. For each $s \in M$, let $A_s = \{a \in A: a \leq s\}$, $P_s = \{a \in A: a \text{ and } s \text{ have no common upper bound}\}$, and $\underline{A_s} = A_s \cup P_s \cup \{a_0\}$. Then we define \[cM = \{A_s:s \in M\} \cup \{\underline{A_s} : s\in M\},\] ordered by inclusion. Additionally we define the \emph{centralization} of $M$ as a subposet of $cM$: \[\underline{M} = \{\underline{A_s}: s \in M\}.\]
\end{definition}

It takes some work to show that $cM$ is in fact a geometric lattice. We will work through the `closure' approach. For a geometric semilattice $M$ with atoms $A$ and a subset $S \subseteq A$, let $J(S)$ be the subposet of $M$ consisting of joins $\bigvee T$, where $T \subseteq S$. It is not hard to check that $J(S)$ is also a geometric semilattice whose rank function agrees with that of $M$.

\begin{lemma} 
\label{closure operator on atoms}
    Define $\cl: 2^A \rightarrow cM$ (where $2^A$ denotes the powerset of $A$) by 

    \[
    \cl(S)=\begin{cases}
        A_{\bigvee S} & a_0 \not \in S, \text{ } \bigvee S \text{ exists} \\
        \underline{A_t}, \text{ } t \text{ maximal in } J(S - a_0) & \text{otherwise}\\
    \end{cases}.
    \]

    Then $\cl$ is a closure operator  satisfying
    \begin{enumerate}[(a)]
        \item $X \subseteq \cl(X).$
        \item If $X \subseteq Y$ then $\cl(X) \subseteq \cl(Y)$
        \item For $X \in cM,$ $\cl(X)=X$.
    \end{enumerate}
\end{lemma}

\begin{proof} 
We first show that $\cl$ is well defined, for which we must show that if $s,t$ are maximal in $J(S-a_0)$, then $\underline{A_s} = \underline{A_t}$. It is a standard fact \cite[Proposition 2.4]{Wachs1985} that $s,t$ have the same rank. Now suppose that we have $a \notin \underline{A_s}$. This implies that $a$ and $s$ have an upper bound which is not $s$, since $a$ is in neither $P_s$ nor $A_s$. Therefore $s < s \vee a$. By the maximality of $s$, we must not have $a \in S$, and since $\operatorname{rank}(s \vee a) > \operatorname{rank}(t)$ we now get that $t < t \vee a$ by property (b) of Definition \ref{def:geometric semilattice}. But this shows that $a \notin \underline{A_t}$. By symmetry, the two sets are equal.

We now prove (a). When $\bigvee X$ exists, certainly $X \subseteq A_{\bigvee X}$. Otherwise, we have $\cl(X) = \underline{A_t}$ for some $t \in M$. By definition we have $a_0 \in \underline{A_t}$, and for any $a \in X - a_0$ the maximality of $t$ in $J(X-a_0)$ implies that either $t = t \vee a$ or $t \vee a$ does not exist. So $a \in A_t$ or $P_t$, and in either case $a \in \underline{A_t}$. 

For (b), let $X,$ $Y \in 2^A$. For some $s \in M$, we have $\cl(X) = A_s$ or $\cl(X) = \underline{A_s}$. Since in the second case $\cl(X)$ being in $\underline{M}$ implies that $\cl(Y)$ is as well, it suffices to show that for some $t \in M$, we have $s \leq t$ and $\cl(Y) = A_t$ or $\cl(Y) = \underline{A_t}$. But since $J(X) \subseteq J(Y)$ there is certainly a maximal element $t$ of $J(Y)$ such that $s \leq t$ as desired.

Finally, for (c) let $X \in cM$. If $X = A_s$, then we have $\bigvee X = s$ so $\cl(X) = A_{\bigvee X} = A_s = X$. If $X = \underline{A_s}$, then by definition for every $a \in X - a_0$, $s \vee a$ does not exist or $a \leq s$. This shows that $s$ is maximal in $J(X - a_0)$, so $\cl(X) = \underline{A_s} = X$.
\end{proof}

\begin{corollary}
   The cone $cM$ is a lattice with join $S \vee T=\cl(S \cup T)$ and meet $S \wedge T=S \cap T.$ 
\end{corollary}
\begin{proof}
Let $S,$ $T$ be elements of $cM$. By Lemma \ref{closure operator on atoms} (a), $S,T \leq \cl(S \cup T).$ Now suppose $Q$ is an element of $cM$ such that $S,T \leq Q$. Therefore $S \cup T \subseteq Q$ and by Lemma \ref{closure operator on atoms} (b), $\cl(S \cup T) \subseteq \cl(Q)$. By Lemma \ref{closure operator on atoms} (c), $\cl(Q)=Q$ so we have $\cl(S \cup T) \leq Q$. Therefore $\cl(S \cup T)$ is the least upper bound of $S$ and $T$.

On the other hand, if $Q \leq S,T$ then certainly $Q \subseteq S \cap T$, so it suffices to show that $S \cap T \in cM$. By Lemma \ref{closure operator on atoms} (b), we have $\cl(S \cap T) \subseteq \cl(S) = S$ and similarly $\cl(S \cap T) \subseteq T$. So $\cl(S \cap T) \subseteq S \cap T$ and we have equality by Lemma \ref{closure operator on atoms} (a). Since the codomain of $\cl$ is $cM$, we have $S \cap T \in cM$ as desired.
\end{proof}

If in some poset we have $s < t$ and $s \leq r \leq t$ implies $r = s$ or $r = t$ (that is, there is nothing in between $s$ and $t$) we say that $t$ \emph{covers} $s$ and write $s \lessdot t$. In order to show that $cM$ is a geometric lattice we will want to understand its covering relations, which the following lemma relates to the closure operator.

\begin{lemma}\label{closure-covering}
    For $S \in cM$ and $a \not \in S,$ $S \lessdot \cl(S \cup \{a\}).$
\end{lemma}
\begin{proof}
We consider three possible cases. If $S = A_s$ and $\cl(S + a) = A_t$ for some $s,t \in M$, then we have $t = \bigvee (S + a) = (\bigvee S) \vee a = s \vee a$. By the semimodularity of $M_{s \vee a}$ we find $\operatorname{rank}(s \vee a) \leq \operatorname{rank}(s) + 1$, so in fact $s \lessdot s \join a$. By the injective order-preserving map $s \mapsto A_s$ we have $A_s \lessdot A_{s \vee a}$ as desired. 

In the next case, if $S = A_s$ and $\cl(S + a) = \underline{A_t}$ we must have that $s$ is maximal in $J(S + a - a_0)$ and in fact $\cl(S + a) = \underline{A_s}$. To see that $A_s \lessdot \underline{A_s}$, suppose that $A_s \leq T \leq \underline{A_s}$ for some $T \in cM$. Since $A_s$ is a maximal subset of $\underline{A_s}$ with a join, either $T = A_s$ or $T = \underline{A_t}$ for some $t \in M$. In the second case, $A_s + a_0 \subseteq \underline{A_t}$ so $\underline{A_s} = \cl(A_s + a_0) \subseteq \cl(\underline{A_t} + a_0) = \underline{A_t}$. Therefore $T = \underline{A_s}$. In either case, $A_s \lessdot \underline{A_s}$.

Finally, we may have $T = \underline{A_s}$ and $\cl(X + a) = \underline{A_t}$. Since we have only added a single atom, by semimodularity of $M_t$ we can choose $s,t$ such that $s \leq t$ and $\operatorname{rank}(s) \geq \operatorname{rank}(t) - 1$. Using $a \notin \underline{A_s}$ we see that this is an equality and in fact $s \lessdot t$. Now suppose there is some in between element $\underline{A_s} \leq \underline{A_r} \leq \underline{A_t}$. If these inequalities are strict, the same argument shows we can pick representatives so that $s \lessdot r \lessdot t$, but this contradicts $s \lessdot t$. Therefore one of the inequalities is actually equality, and we have $\underline{A_s} \lessdot \underline{A_t}$.
\end{proof}

\begin{proposition}
   The cone $cM$ is a geometric lattice.
\end{proposition}
\begin{proof}
We have already established that $cM$ is a lattice. As a sublattice of the boolean lattice on $A + a_0$ including every singleton as an atom, $cM$ is atomistic. To show that $cM$ is ranked and semimodular, it suffices to show that whenever $S \wedge T \lessdot S,T$ in $cM$, we have $S,T \lessdot S \vee T$ \cite[Proposition 3.3.2]{Stanley_2012}. If this is the case, let $a \in S - (S \wedge T)$. By the properties of $\cl$ and the assumption $S \wedge T \lessdot S$, we have $\cl(S \wedge T + a) = S$. Since $S \wedge T \subseteq T$, we have $S = \cl(S \wedge T + a) \subseteq \cl(T + a)$. Thus $\cl(T + a)$ is an upper bound on both $S$ and $T$, and by Lemma \ref{closure-covering}, $T \lessdot \cl(T + a)$. Symmetrically, for $b \in T - (S \wedge T)$, $\cl(S + b)$ is an upper bound on $S$ and $T$ and $S \lessdot \cl(S +b)$. It follows that $\cl(S + b) = \cl(T + a) = S \join T$, and thus $S,T \lessdot S \join T$.
\end{proof}

We should now confirm that this combinatorial cone corresponds to the geometric one. The isomorphism between the two is fairly obvious but somewhat tedious to verify, so we will rely on a uniqueness result from \cite{Wachs1985}.

\begin{proposition} \label{equivalence}
    Let $\Aa$ be an arrangement in $\Rn$. Then there is a poset isomorphism $\varphi: \Ll(c\Aa) \to c\Ll(\Aa)$ such that $\varphi(H_0) = \{a_0\}$. 
\end{proposition}
\begin{proof}
In fact the isomorphism also sends each $cH_i$ to $\{H_i\}$ and otherwise is determined by atomicity. Theorem 3.2 of \cite{Wachs1985} states that any geometric semilattice $M$ is embedded as $L - L^a$ for some unique geometric lattice $L$ and atom $a$ of $L$. Lemma \ref{lem:cone_geo} shows that $\Ll(\Aa) = \Ll(c\Aa) - \Ll(c\Aa)^{H_0}$, where $\Ll(c\Aa)$ is a geometric lattice and $H_0$ an atom. On the other hand, the map $X \mapsto A_X$ gives an isomorphism $\Ll(\Aa) \to c\Ll(\Aa) - c\Ll(\Aa)^{\{a_0\}}$, where $c\Ll(\Aa)$ is a geometric lattice with atom $\{a_0\}$. By the uniqueness of the geometric lattice for which this embedding holds, these pairs $(\Ll(c\Aa), H_0)$ and $(c\Ll(\Aa), \{a_o\})$ must be isomorphic.
\end{proof}

With this in hand we can complete the proof of Theorem \ref{degrees of freedom theorem}, since by its definition $c\Ll(\Aa)$ is clearly a combinatorial invariant.

\begin{theorem}\label{combinatorial_inv}
    The poset of the centralization $\Ll(\underline{\Aa})$ is a combinatorial invariant of $\Aa$. Moreover for any $V \in \Ll(\underline{\Aa})$, $\Ll(\underline{\Aa}^V)$ and $\Ll(\Aa_V)$ are combinatorial invariants as well.
\end{theorem}
\begin{proof}
By definition, $c\Ll(\Aa)$ is a combinatorial invariant, and the isomorphism of Proposition \ref{equivalence} shows that $\Ll(c\Aa)$ is as well. Moreover, the atom $H_0$ of $\Ll(c\Aa)$ is $\varphi^{-1}(a_0)$, so we can identify it in $\Ll(c\Aa)$. Lemma \ref{lem:cone_geo} lets us construct $\Ll(\underline{\Aa})$ and $\Ll(\Aa_V)$ from $\Ll(c\Aa)$ through the combinatorial operations of taking subposets, order filters, and order ideals. Finally we have $\Ll(\underline{\Aa}^V) \simeq \Ll(\underline{\Aa})^V$, which we can obtain from $\Ll(\underline{\Aa})$.
\end{proof}

\section{The level distribution for semilattices} \label{sec:semilattice}

We now observe that the construction of $cM$ allows us to generalize the definition of $r_\ell$ to all geometric semilattices, including those that are not realizable as the intersection poset of a real hyperplane arrangement. We make a series of definitions in analogy to the components of Theorem \ref{degrees of freedom theorem}.

\begin{definition}
    For a geometric semilattice $M$, let $cM$ and $\underline{M}$ be as in Definition \ref{combo_cone}, and for $S \in \underline{M}$ let $M_S = cM_S - \underline{M}$ be the \emph{localization} of $M$ to $S$.
\end{definition}

To define $r(M)$ and $b(M)$ we will use the characteristic polynomial of $M$, which we define to match the characteristic polynomial for a hyperplane arrangement.
\begin{definition}
    For a geometric semilattice $M$ of rank $n$, define the \emph{characteristic polynomial of $M$} to be 
    \[\chi_M(t) = \sum\limits_{s \in M} \mu(\hat 0, s)t^{n-\operatorname{rank}(s)}.\]
    Also define the values
     \begin{equation}
        r(M) = (-1)^n\chi_M(-1), \hspace{10pt} b(M) = (-1)^{n}\chi_M(1).
    \end{equation}
\end{definition}

Using the definitions so far, we can define the numbers $r_\ell(M)$ for any geometric semilattice.

\begin{definition}
    For a geometric semilattice $M$ of rank $n$, for $0 \leq \ell \leq n$ we define the values:
    \begin{equation}
        r_\ell(M) = \sum\limits_{S \in \underline{M}_{n - \ell}} r(\underline{M}^S) b(M_S).
    \end{equation}
\end{definition}

As a consequence of the fact that $\mu(\hat 0, s)$ alternates in sign with the rank of $s$, it is not hard to determine that all these values are positive, and the definition also causes $r_0(M) = b(M)$. As a geometric fact for real hyperplane arrangements we would hope that the sum of the $r_\ell(M)$ is $r(M)$ in general. We can obtain this fact from Theorem \ref{chi_semilattice} by evaluating at $t = -1$ and using the definitions, taking some care with signs. These properties suggest that $r(M)$ and $r_\ell(M)$ enumerate some objects associated to $M$ analogous to the regions of an arrangement, and that these objects can be naturally partitioned according to their `level'. We are not aware of such interpretations for non-realizable geometric semilattices, presenting an avenue for future exploration. Having made the appropriate definitions, it is now not too difficult to prove Theorem \ref{chi_semilattice}.

\begin{proof}[Proof of Theorem \ref{chi_semilattice}]
Recall that a (defining) property of the M\"obius function is that for a ranked poset $L$ and $S,T \in L$ we have $\sum_{S \leq U \leq T} \mu(U,T) = \delta_{ST}$. Now identify $M$ with $cM - \underline{M}$, and for an element $U \in M$, let $\underline{U}$ be the unique element of $\underline{M}$ such that $U \lessdot \underline{U}$. We expand the right hand side of Theorem \ref{chi_semilattice} as
\begin{align*}
    p(t) = \sum\limits_{S \in \underline{M}} \chi_{\underline{M}^S}(t)\chi_{M_S}(1) &= \sum \limits_{S \in \underline{M}} \left[\left(\sum\limits_{T \in \underline{M}^S} t^{(n+1) - \operatorname{rank}(T)}\mu_{cM}(S,T)\right)\left(\sum\limits_{U \in M_S} \mu_{cM}(\hat 0, U)\right)\right].\\
\end{align*}
We now reorder as a nested sum for $U$, $T$, and then $S$ such that $U \lessdot \underline{U} \leq S \leq T$ to obtain
\begin{align*}
    p(t) &= \sum\limits_{U \in M} \left[\mu_{cM}(\hat 0 , U) \sum\limits_{\underline{U} \leq T \in \underline{M}} \left(t^{(n+1)-\operatorname{rank}(T)} \sum\limits_{\underline{U} \leq S \leq T} \mu_{cM}(S,T)\right)\right]. 
\end{align*}

The innermost sum collapses to $\delta_{\underline{U}T}$, so the only remaining term of the middle sum is the term for $\underline{U}$. Using the fact that $\text{rk}(\underline{U}) = \text{rk}(U) + 1$, and $\mu_{cM}(\hat 0,U) = \mu_M(\hat 0, U)$ we have
\begin{align*}
    p(t) = \sum\limits_{U \in M}\mu_M(\hat 0,U) t^{n - \text{rk}(U)} = \chi_M(t),
\end{align*}
as we set out to prove.
\end{proof}

In the case of real hyperplane arrangements, we can use the bijection of Theorem \ref{bijection} to give a sum over regions.

\begin{theorem} \label{chi_thm}
    For a hyperplane arrangement $\Aa$ in $\RR^n$, we have 
    \begin{equation}
        \chi_\Aa(t) = (-1)^n\sum\limits_{R \text{ a region of $\Aa$}} \chi_{\underline{\Aa}^{V_R}}(t)/\chi_{\underline{\Aa}^{V_R}}(-1).
    \end{equation}
\end{theorem}
\begin{proof}
When the geometric semilattice in Theorem \ref{chi_semilattice} is the intersection poset of a hyperplane arrangement $\Aa$, we have
\[\chi_\Aa(t) = \sum\limits_{V \in \Ll(\underline{\Aa})} \chi_{\underline{\Aa}^V}(t)\chi_{A_V}(1) = (-1)^n\sum\limits_{V \in \Ll(\underline{\Aa})}r(\underline{\Aa}^V)b(\Aa_V)(\chi_{\underline{\Aa}^V}(t)/\chi_{\underline{\Aa}^V}(-1)),\]
substituting by Zaslavsky's Theorem (\ref{Zaslavsky}). By Theorem \ref{bijection}, the factor $r(\underline{\Aa}^V)b(\Aa_V)$ counts the number of regions $R$ such that $V_R = V$. Therefore we can break each term for $V$ into a sum over its regions, completing the proof.
\end{proof}

In another application, for certain semilattices Theorem \ref{chi_semilattice} has an especially nice form.

\begin{definition} \label{uniform}
    Let $L$ be a geometric lattice of rank $n$. If there is a sequence of geometric lattices $L_0,\ldots,L_n$ such that for any $S \in L_{n-d}$ we have $L^S \simeq L_d$, we call $L$ \emph{restriction-uniform} and define for $0 \leq d \leq n$ the polynomial $\tilde\chi_L^d(t) = \chi_{L_d}(t)/\chi_{L_d}(-1)$ . 
\end{definition}

Examples of restriction-uniform geometric lattices include the partition and signed partition lattices, the boolean lattice, and its $q$-analogs. When a semilattice $M$ has restriction-uniform centralization, we can gather the terms in Theorem \ref{chi_semilattice} by the dimension of $S$ and rearrange to obtain the following corollary.

\begin{corollary}\label{cor:uniform}
    Let $M$ be a geometric semilattice such that $\underline{M}$ is restriction-uniform. Then we have 
    \begin{equation}
        \chi_M(t) = \sum\limits_{\ell = 0}^n r_\ell(M)\tilde\chi_{\underline{M}}^\ell(t).
    \end{equation}
\end{corollary}

This corollary motivated our discovery of Theorem \ref{chi_semilattice}, as the special case where $\underline{M}$ is $\Pi_n$ captures \emph{nondegenerate deformations} of the braid arrangement, recovering the main result of \cite{zhang_2024} and generalizing a result from \cite{chen_2026} as we will see in the next section.

\section{Deformations of the Braid Arrangement} \label{sec:exp}
In this section we apply Theorem \ref{degrees of freedom theorem} to the class of exponential arrangements, which we will define shortly. Zaslavsky gave the specialization of his equation to such arrangements in \cite{Zaslavsky_2003}, in particular giving an alternate enumeration for the extended Shi arrangement, which Ehrenborg studied when first defining level. More recently Chen, Fu, Wang, and Yang produced generating function identities regarding a very general subclass of exponential arrangements \cite{chen_2026}. We show that Zasklavsky's specialization implies these results for all exponential arrangements. Note that terminology differs in some sources, Zaslavsky calls deformations of the braid arrangement \emph{affinographic} and refers to exponential sequences of arrangements as possessing a ``kind of uniformity". Chen, Fu, Wang, and Yang also gave an expression for $\chi_\Aa(t)$ in terms of the $r_\ell(\Aa)$ which can be seen as a generalization of Ehrenborg's original result, and which was in turn generalized by Zhang \cite{zhang_2024}. We show that Corollary \ref{uniform} can be applied to give a new proof of this.

Exponential arrangements are closely related to the \emph{braid} (or Coexter type $A$) arrangement in $\Rn$, given by 
\[
    \text{Braid}(n) = \{x_i - x_j = 0 : 1 \leq i < j \leq n\}.
\]
Let us review some useful facts about the braid arrangement. The braid arrangement is central, so all of its regions have $\ell(R) = n$. There is a simple bijection between permutations $\omega \in S_n$ and the regions of Braid$(n)$ given by
\[\omega \in S_n \leftrightarrow \{(x_1,\ldots,x_n) \in \Rn : x_{\omega(1)} > x_{\omega(2)} > \cdots > x_{\omega(n)}\},  \]
and therefore $r_n(\text{Braid}(n)) = r(\text{Braid}(n)) = n!$. For $n > 0$, the braid arrangement is nonessential with its smallest flat being the 1-dimensional subspace $x_1 = x_2 = \cdots = x_n$. The intersection poset of the braid arrangement is isomorphic to the lattice $\Pi_n$ of set partitions on $\{1,\ldots,n\}$, ordered by refinement. We can interpret a set partition $B = \{B_1,\ldots,B_k\}$ as a flat of dimension $k$ where we have $x_i = x_j$ exactly when $i$ and $j$ are in the same part of $B$. This lattice has a self-similarity with its principal order filters such that  $(\Pi_n)^B \simeq \Pi_{|B|}$. Therefore we have $\Ll(\text{Braid}(n)^F) \simeq \Ll(\text{Braid}(\dim(F)))$ for any flat $F$. Principal order ideals of $\Pi_n$ also have the nice property that they are a product of smaller copies: $(\Pi_n)_B \simeq \prod_{i=1}^{|B|}\Pi_{|B_i|}$. Therefore the intersection poset of a localization of $\text{Braid}(n)$ is a product of intersection posets of smaller braid arrangements.

Exponential sequences of arrangements are defined in such a way that the properties we have just highlighted are in some sense preserved. Call an arrangement $\Aa$ a \emph{nondegenerate deformation} of a central arrangement $\Bb$ if $\underline{\Aa} = \Bb$, that is each of the hyperplanes of $\Aa$ is parallel to a hyperplane of $\Bb$, and for each hyperplane of $\Bb$ there is a hyperplane of $\Aa$ parallel to it.

\begin{definition} \label{def:exp}
    An exponential sequence of arrangements $\Aa_1,\Aa_2,\ldots$ is a sequence such that 
    \begin{enumerate}
        \item $\Aa_n$ is a nondegenerate deformation of $\text{Braid}(n)$.
        \item For any subset $S \subseteq \{1,\ldots,n\}$, let $V_S = \{x \in \Rn : \text{for all } i,j \in S, x_i = x_j\}$. Then \[\Ll((\Aa_n)_{V_S}) \simeq \Ll(\Aa_{|S|}).\]
    \end{enumerate}
\end{definition}
\begin{figure}
    \centering
    \includegraphics[width=4cm]{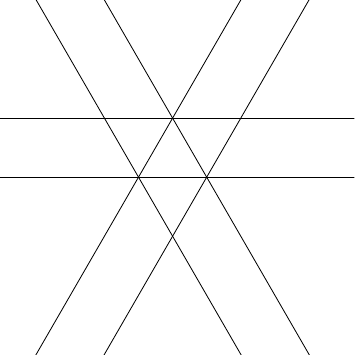}
    \caption{$\text{Shi}(3) = \{x_i - x_j = 0,1 : 1 \leq i < j \leq 3\}$}
    \label{Shi(3)}
\end{figure}
Note that the braid arrangement itself is an exponential arrangement. Other popular examples are the Catalan, Semiorder, and Shi arrangements (see Fig. \ref{Shi(3)}). The first property in this definition ensures that $\underline{\Aa_n} \simeq \text{Braid}(n)$, and therefore that $\Aa_n$ is nonessential with $\dim(\operatorname{ess}(\Aa_n)) = n - 1$ so that $b(\Aa_n) = r_1(\Aa_n)$. The second property causes the principal order ideals of $\Ll(\Aa_n)$ to decompose similarly to those of the braid arrangement.
\begin{lemma}
    Let $\Aa_n$ be an exponential sequence of arrangements. If $F$ is a flat of $\Ll(\underline{\Aa_n}) \simeq \Pi_n$ with associated partition $B$, then we have $\Ll((\Aa_n)_F) \simeq \prod_{i=1}^{|B|} \Ll(\Aa_{|B_i|})$.
\end{lemma}
\begin{proof}
    The flat $F$ consists of all points where $x_i = x_j$ when $i \sim_B j$. Therefore $(\Aa_n)_F$ consists of all hyperplanes of the form $x_i - x_j = a$ for $1 \leq i < j \leq n$ and $i \sim_B j$. We may associate each hyperplane of $(\Aa_n)_F$ with the part of $B$ such that the indices of the coordinates in its definition both appear in that part, and the set of hyperplanes associated to the part $B_r$ is $(\Aa_n)_{V_{B_r}}$ where $V_{B_r}$ is the subspace defined in Definition \ref{def:exp}. For two hyperplanes $H_1,H_2$ associated to different parts of $B$, we have $H_1 \perp H_2$. It follows that \[(\Aa_n)_F = (\Aa_n)_{B_1} \times (\Aa_n)_{B_2} \times \cdots \times (\Aa_n)_{B_k}\] is a product arrangement in the sense of \cite[Definition 2.13]{orlik_terao_1992}. By \cite[Proposition 2.14]{orlik_terao_1992} the intersection poset is a product of posets, and by the second property of exponential arrangements each factor is isomorphic to $\Ll(\Aa_{|B_i|}).$
\end{proof}

We can simplify Theorem \ref{degrees of freedom theorem} in the case of exponential sequences of arrangements by applying the properties we have just collected. The following is equivalent to the $x^0$ term of Equation 6.2 in \cite{Zaslavsky_2003}.

\begin{theorem} \label{thm:exp}
    For $\Aa_n$ an exponential sequence of arrangements, 
    \begin{equation}
        r_\ell(\Aa_n) = \ell! \sum\limits_{B \in (\Pi_n)_{n - \ell}} \prod\limits_{i=1}^\ell b(\Aa_{|B_i|}).
    \end{equation}
\end{theorem}
\begin{proof}
    This follows in a straightforward manner from Theorem \ref{degrees of freedom theorem}, making substitutions according to the observations above. Additionally we use that fact that the characteristic polynomial of a product poset is the product of characteristic polynomials, so $b(\prod \Aa_i) = \prod b(\Aa_i)$.
\end{proof}

Theorem \ref{thm:exp} can be restated in terms of generating functions, generalizing Theorem 1.1  of \cite{chen_2026}. For an exponential sequence of arrangements $\Aa$, let \[R_\ell(\Aa;t) = \sum\limits_{n \geq 0} r_\ell(\Aa_n)\frac{x^n}{n!}.\] Then $R_\ell(\Aa, t)$ exhibits a binomial type identity.
\begin{theorem} \label{thm:binom}
    For any exponential sequence of arrangements $\Aa$ and $\ell \geq 0$, we have
    \begin{equation}
        R_\ell(\Aa;t) = R_1(\Aa,t)^\ell.
    \end{equation}
    Moreover,  for any fixed integers $\ell_1,\ell_2 \geq 0$,
    \begin{equation}
    r_{\ell_1 + \ell_2}(\Aa_n) = \sum\limits_{i=0}^n \binom{n}{i}r_{\ell_1}(\Aa_i)r_{\ell_2}(\Aa_{n-i})
    \end{equation} holds for any $n \geq 0$.
\end{theorem}
\begin{proof}
    After some manipulation of exponential generating functions, it suffices for the both equalities to show that 
    \begin{equation} \label{eqn:binom}
        \frac{1}{n!}r_\ell(\Aa_n) = \sum\limits_{\substack{n_1  + \cdots + n_\ell= n \\ n_i > 0}} \prod\limits_{i=1}^\ell\frac{r_1(\Aa_{n_i})}{n_i!}.
    \end{equation}
    This can be obtained from Theorem \ref{thm:exp} by re-grouping terms of the sum with identical summands, so that the index set is over sums $n_1 + \cdots + n_\ell = n$ rather than partitions with $\ell$ parts. Since the $\binom{n}{n_1,\ldots,n_\ell}$ partitions with part sizes $n_1,\ldots,n_\ell$ in some order correspond to $\ell!$ such sums, this incurs a multiple of $\frac{1}{\ell!}\binom{n}{n_1,\ldots,n_\ell}$. This cancels the $\frac{1}{n!}$ on the left of Equation \eqref{eqn:binom} and the $\ell!$ in Theorem \ref{thm:exp}, and we associate each $\frac{1}{n_i!}$ with a term of the product to obtain the right had side.
\end{proof}

In \cite{chen_2026} the authors derive an additional consequence of this Theorem by making use of a theorem of Stanley which relates the exponential generating functions for the characteristic polynomial and number of regions of an exponential sequence of arrangements \cite[Theorem 5.17]{stanley2007introduction}. This allows them to show that the values $r_\ell(\Aa)$ serve as coefficients for the characteristic polynomial in the binomial basis. Theorem \ref{thm:binom} means that the same proof applies to all exponential sequences of arrangements. In fact the theorem holds for all nondegenerate deformations of the braid arrangement which Zhang proves in \cite{zhang_2024} using deletion-restrictions techniques. We can give another proof of this through Corollary \ref{cor:uniform}.

\begin{theorem}{\cite[Theorem 1.1]{zhang_2024}} \label{thm:chi_exp}
    Let $\Aa$ be a hyperplane arrangement such that $\underline{\Aa} \simeq \text{Braid}(n)$. Then we have
\begin{equation} \label{eqn:chi_exp}
    \chi_{\Aa}(t) = \sum\limits_{\ell = 0}^n (-1)^{n - \ell}r_\ell(\Aa)\binom{t}{\ell}.
\end{equation}
    Conversely, the number of regions of level $\ell$ in $\Aa$ can be expressed as
\begin{equation} \label{coeffs}
    r_\ell(\Aa) = (-1)^{n}\sum\limits_{k=0}^\ell (-1)^k \binom{\ell}{k} \chi_{\Aa}(k) = (-1)^n\Delta^\ell \chi_{\Aa}(t)|_{t = 0},
\end{equation}
where $\Delta$ is the difference operator $p(t) \mapsto p(t) - p(t+1)$.
\end{theorem}
\begin{proof}
    As noted earlier, for any flat  $V \in \Ll(\text{Braid}(n))$ of rank $n-\ell$ and therefore dimension $\ell$, we have $\Ll(\text{Braid}(n))^V \simeq \Ll(\text{Braid}(\ell))$. Therefore the braid arrangement is restriction-uniform (Definition \ref{uniform}). 
    Since $\chi_{\text{Braid}(\ell)}(t) = t(t-1)\cdots(t-\ell+1)$, we have $\tilde\chi_{\underline{\Aa}}^{n-\ell}(t) = (-1)^\ell\binom{t}{\ell}$. Therefore Equation \eqref{eqn:chi_exp} is just this case of Corollary \ref{cor:uniform}. Because the (alternating) binomial coefficients are a basis for the space of polynomials, we can solve for each coefficient $r_\ell(\Aa)$, and this is the first equality of Equation \eqref{coeffs}. One can check the second equality simply by expanding the right hand side.
\end{proof}

This shows that in the case of exponential arrangements, $r_\ell$ actually is a characteristic invariant. We give the second expression for $r_\ell(\Aa)$ in Equation \eqref{coeffs} because it can be seen as a generalization of the main result of \cite{ehrenborg_level}, once that theorem is restricted to count regions, rather than faces of any dimension. In this way Ehrenborg's result is unified with the other results of this kind. We note that this also gives an alternative proof of the fact that the Shi and Ish arrangements share values of $r_\ell$ as studied in \cite{Armstrong_Rhoades_2011}, since both are nondegenerate deformations of the braid arrangement and that paper applies the finite fields method of \cite{athanasiadis_1996} to show that they have the same characteristic polynomial.
Zhang also gives a similar expression for nondegenerate deformations of the Coxeter arrangement of type $B$, which can also be seen as a special case of Corollary \ref{cor:uniform}.

\printbibliography[category=cited]

@article{Armstrong_Rhoades_2011, title={The Shi Arrangement and the Ish Arrangement}, volume={DMTCS Proceedings vol. AO,...}, DOI={10.46298/dmtcs.2890}, number={Proceedings}, journal={Discrete Mathematics \&amp; Theoretical Computer Science}, author={Armstrong, Drew and Rhoades, Brendon}, year={2011}, month={1}}

@book{zaslavsky_1975,
  title={Facing up to arrangements: Face-count formulas for partitions of space by hyperplanes},
  author={Zaslavsky, Thomas},
  volume={154},
  year={1975},
  publisher={American Mathematical Soc.}
}

@article{athanasiadis_1996,
title = {Characteristic Polynomials of Subspace Arrangements and Finite Fields},
journal = {Advances in Mathematics},
volume = {122},
number = {2},
pages = {193-233},
year = {1996},
issn = {0001-8708},
doi = {https://doi.org/10.1006/aima.1996.0059},
url = {https://www.sciencedirect.com/science/article/pii/S0001870896900596},
author = {Christos A. Athanasiadis}
}

@inproceedings{zhang_2024,
  title={Characteristic Polynomials of Deformations of Coxeter Arrangements via levels of regions},
  author={Ningxin Zhang},
  year={2024},
  url={https://api.semanticscholar.org/CorpusID:273850458}
}

@article{zaslavsky_1977,
title = {A combinatorial analysis of topological dissections},
journal = {Advances in Mathematics},
volume = {25},
number = {3},
pages = {267-285},
year = {1977},
issn = {0001-8708},
doi = {https://doi.org/10.1016/0001-8708(77)90076-7},
url = {https://www.sciencedirect.com/science/article/pii/0001870877900767},
author = {Thomas Zaslavsky}
}

@article{stanley2007introduction,
  title={An introduction to hyperplane arrangements},
  author={Stanley, Richard P and others},
  journal={Geometric combinatorics},
  volume={13},
  pages={389--496},
  year={2007}
}

@misc{ardila2004,
      title={Computing the Tutte polynomial of a hyperplane arrangement}, 
      author={Federico Ardila},
      year={2004},
      eprint={math/0409211},
      archivePrefix={arXiv},
      primaryClass={math.CO},
      url={https://arxiv.org/abs/math/0409211}, 
}

@incollection{ardila2022tutte,
  title={Tutte polynomials of hyperplane arrangements and the finite field method},
  author={Ardila, Federico},
  booktitle={Handbook of the Tutte Polynomial and Related Topics},
  pages={532--550},
  year={2022},
  publisher={Chapman and Hall/CRC}
}

@misc{chen2025level,
      title={Regions of Level $\ell$ of Catalan/Semiorder-Type Arrangements}, 
      author={Yanru Chen and Suijie Wang and Jinxing Yang and Chengdong Zhao},
      year={2025},
      eprint={2410.10198},
      archivePrefix={arXiv},
      primaryClass={math.CO},
      url={https://arxiv.org/abs/2410.10198}, 
}

@article{chen_2026,
title = {Level of regions for deformed braid arrangements},
journal = {Journal of Combinatorial Theory, Series A},
volume = {217},
pages = {106077},
year = {2026},
issn = {0097-3165},
doi = {https://doi.org/10.1016/j.jcta.2025.106077},
url = {https://www.sciencedirect.com/science/article/pii/S009731652500072X},
author = {Yanru Chen and Houshan Fu and Suijie Wang and Jinxing Yang}
}

@Article{Wachs1985,
author={Wachs, Michelle L.
and Walker, James W.},
title={On geometric semilattices},
journal={Order},
year={1985},
month={12},
day={01},
volume={2},
number={4},
pages={367-385},
issn={1572-9273},
doi={10.1007/BF00367425},
url={https://doi.org/10.1007/BF00367425}
}

@article{ehrenborg_level,
title = {Counting faces in the extended Shi arrangement},
journal = {Advances in Applied Mathematics},
volume = {109},
pages = {55-64},
year = {2019},
issn = {0196-8858},
doi = {https://doi.org/10.1016/j.aam.2019.05.002},
url = {https://www.sciencedirect.com/science/article/pii/S0196885819300855},
author = {Richard Ehrenborg}
}

@misc{ardila2004_semimat,
      title={Semimatroids and their Tutte polynomials}, 
      author={Federico Ardila},
      year={2004},
      eprint={math/0409003},
      archivePrefix={arXiv},
      primaryClass={math.CO},
      url={https://arxiv.org/abs/math/0409003}, 
}

@article{Postnikov_2000,
title = {Deformations of Coxeter Hyperplane Arrangements},
journal = {Journal of Combinatorial Theory, Series A},
volume = {91},
number = {1},
pages = {544-597},
year = {2000},
issn = {0097-3165},
doi = {https://doi.org/10.1006/jcta.2000.3106},
url = {https://www.sciencedirect.com/science/article/pii/S0097316500931066},
author = {Alexander Postnikov and Richard P. Stanley}
}

@misc{hetyei2025,
      title={Labeling regions in deformations of graphical arrangements}, 
      author={Gábor Hetyei},
      year={2025},
      eprint={2312.06513},
      archivePrefix={arXiv},
      primaryClass={math.CO},
      url={https://arxiv.org/abs/2312.06513}, 
}

@phdthesis{Dorpalen-Barry_2021, title={Cones of Hyperplane Arrangments}, school={University of Minnesota}, url={https://hdl.handle.net/11299/224994}, author={Dorpalen-Barry, Galen Anna}, year={2021}}

@book{orlik_terao_1992,
author = {Orlik, Peter. and Terao, Hiroaki.},
address = {Berlin, Heidelberg},
booktitle = {Arrangements of Hyperplanes},
isbn = {9783662027721},
publisher = {Springer Berlin Heidelberg},
language = {eng},
series = {Grundlehren der mathematischen Wissenschaften, A Series of Comprehensive Studies in Mathematics, 300},
title = {Arrangements of Hyperplanes },
url = {https://openurl.ac.uk/ukfed:bat?u.ignore_date_coverage=true&rft.mms_id=991003560614302761},
year = {1992},
}

@book{Stanley_2012, place={New York}, title={Enumerative combinatorics. volume 1}, publisher={Cambridge University Press}, author={Stanley, Richard P.}, year={2012}}

@article{stanley_1995,
 ISSN = {00278424},
 URL = {http://www.jstor.org/stable/38724},
 author = {Richard P. Stanley},
 journal = {Proceedings of the National Academy of Sciences of the United States of America},
 number = {6},
 pages = {2620--2625},
 publisher = {National Academy of Sciences},
 title = {Hyperplane Arrangements, Interval Orders, and Trees},
 urldate = {2025-09-03},
 volume = {93},
 year = {1996}
}

@article{Budur_2022,
title = {On contact loci of hyperplane arrangements},
journal = {Advances in Applied Mathematics},
volume = {132},
pages = {102271},
year = {2022},
issn = {0196-8858},
doi = {https://doi.org/10.1016/j.aam.2021.102271},
url = {https://www.sciencedirect.com/science/article/pii/S0196885821001093},
author = {Nero Budur and Tran Quang Tue}
}

@article{Zaslavsky_2003,
author = {Zaslavsky, T.},
year = {2003},
month = {04},
pages = {63-80},
title = {Faces of a Hyperplane Arrangement Enumerated by Ideal Dimension, with Application to Plane, Plaids, and Shi},
volume = {98},
journal = {Geometriae Dedicata},
doi = {10.1023/A:1024029318990}
}

\end{document}